%% file: main_ieee.tex
\documentclass[lettersize,journal]{IEEEtran}
\usepackage{amsmath,amsfonts, amssymb,amsthm, bm}
\usepackage{algorithmic}
\usepackage{algorithm}
\usepackage{array}
\usepackage{multirow, booktabs}
\usepackage[caption=false,font=normalsize,labelfont=sf,textfont=sf]{subfig}
\usepackage{textcomp}
\usepackage{stfloats}
\usepackage{url}
\usepackage{verbatim}
\usepackage{graphicx}
\usepackage{cite}
\usepackage{enumitem}
\usepackage{newtxmath}
\usepackage{setspace}
\usepackage{booktabs}
\usepackage{csvsimple}
\usepackage{siunitx}
\usepackage{makecell}
\usepackage[export]{adjustbox}
\usepackage{comment}
\usepackage{cases}
\usepackage{hyperref}
\allowdisplaybreaks
\hypersetup{
     colorlinks =true,
     linkcolor = blue,
     filecolor = blue,
     citecolor = magenta,      
     urlcolor = magenta,
     }
\usepackage{setspace}
\usepackage{siunitx}
\usepackage[T1]{fontenc}
\newcommand{\doe}[2]{\frac{\partial #1}{\partial #2}}

\theoremstyle{remark}
\newtheorem{remark}{Remark}
\newtheorem{definition}{Definition}

\begin{document}

\title{Steady Flow of Natural Gas in Pipeline Networks via Solution of a Nonlinear Differential-Algebraic System of Equations}



\author{Shriram Srinivasan and Kaarthik Sundar
\thanks{The authors are staff-scientists at Los Alamos National Laboratory, Los Alamos, NM 87544.
email: \texttt{\{shrirams,kaarthik\}@lanl.gov}}}



\maketitle

\begin{abstract}
\input{sections/abstract}
\end{abstract}


\input{sections/introduction}
\input{sections/pipe}

\input{sections/network}
\input{sections/results}

\input{sections/conclusion}

\bibliographystyle{plain}
\bibliography{refs}
\end{document}

%% file: sections/abstract.tex
In the consideration of steady-state flow of gas in pipeline networks, the exclusion of gravity and nonlinear inertial effects (convective acceleration) leads to a fortuitous simplification in the governing equations  to yield a system of nonlinear  algebraic equations. 
Consequently, there are no studies that quantify the effect of gravity and inertial effects on the flow of gas in pipeline networks or delineate regimes of flow conditions wherein the effects are significant or negligible.
In addressing this need, we consider the steady-state flow equations in pipeline networks  without neglecting the gravitational and inertial terms and in place of a system of algebraic equations (one for each pipe), this approach results in a nonlinear system of first-order ordinary differential equations (ODEs) which are coupled through algebraic equations that appear in the form of boundary conditions on the pressure and  balance of mass flows at either end.
One of our main contributions in this article is to demonstrate how the Newton-Raphson algorithm can still be used to solve the  coupled nonlinear  differential-algebraic system by utilizing  the appropriate forward sensitivity ODEs to evaluate the Jacobian terms arising in the iterative scheme.
We also propose a variable transformation that alleviates the poor scaling of the ODE, and we introduce a two-point collocation scheme as a coarse approximation of the system from which to find initial guesses for the Newton iterations. 
Simulation studies were conducted for a single pipe as well as a large-scale pipeline network with real data. From these studies, we concluded that while the effect of gravity is important, the inertial effect was negligible in all cases. 
The proposed methodology is applicable to a wide class of pipeline and thermal-fluid networks beyond natural gas, including liquid pipelines, hydrogen transport, and district heating systems.

%% file: sections/introduction.tex
\section{Introduction}
\label{sec:introduction}
The task of simulating gas flow in a network of pipes occurs in a variety of situations and applications, such as the optimization of network operations, estimation of delivery capacity, as well as contingency analysis and planning for natural hazards. 
While it is undeniable that the gas flow and pressures vary with time (daily as well as seasonally due to weather patterns), the case of steady or quasi-static flow  is nevertheless important for reasons other than the obvious  mathematical simplification that it introduces  in the governing equations. 
The length and span of gas pipeline networks mean that measurements of pressure and flow are sparse, and often reported as averages over a meaningful time interval, over which variations are often small enough to be neglected. 
Moreover, normal operating conditions engender flows in the subsonic regime which can be approximated as steady, unless the phenomena of interest are events such as the closing of a valve or changes at a control point that occur on small timescales.

Thus, the simulation of steady flow in pipeline networks has drawn much interest from engineers and researchers in the field.  
Despite the simplification introduced by the assumption of steady flow, the mathematical problem remains challenging due to the following reasons: 
\begin{enumerate}[label=(\arabic*)]
    \item The  (\emph{non-ideal}) equation of state for natural gas  that has a nonlinear relationship between pressure and density
    \item Nonlinear inertial effects due to convective acceleration 
    \item Gravitational effects introduced due to changes in elevation by non-horizontal pipelines 
    \item The sheer size  and scale of the pipeline network 
\end{enumerate}

Published research thus far has met varying degrees of success at addressing \emph{at most} two of these four challenges.
Most studies neglect the inertial effect since it was shown to be negligible in comparison to the viscous drag for a single pipe \cite{Osiadacz1984Jan}. 
A host of studies have made further simplifying assumptions in the form of the ideal gas equation of state and horizontal pipelines to focus on theoretical issues \cite{Dvijotham2015Jun,ss-soln-existence}  or propose  computational methods that can handle large networks with limited success \cite{ojha2017solving, Singh2019, Singh2020}. 
Some research efforts \cite{nrsolver,network-partitioning} were successful in simulating steady-state flows in large-scale networks with both ideal and non-ideal equations of state while neglecting inertial and gravitational effects. 
On the other hand, the authors in \cite{Gugat2015}  incorporate the inertial and gravitational effects assuming the ideal gas equation of state but restrict their attention to flow in a single pipe while in \cite{Gugat2018May}, non-ideal gases with inertia are considered in the analysis of flow on networks but gravity is neglected and the focus is not on computational methods for solving the problem in the case of large networks. 
In summary, \emph{there has been no attempt} to simulate steady flow of natural gas in large-scale pipeline networks that also takes into account (1), (2), and (3) above.

Since the difference between the ideal and non-ideal equations of state and their relevance to pipeline network simulations is well-known \cite{GyryaZlotnik2017,nrsolver}, we shall not discuss that further but instead focus on why gravitational and inertial effects \emph{could} be important and ought to be incorporated in simulations.
As pipeline networks often have a wide expanse (sometimes several states in the continental US), it is quite likely that sections of the pipeline would span both mountainous terrain and flat plains so that some of the pipelines must be non-horizontal.
Depending on the inclination and the direction of flow in the pipeline, gravity could either accelerate the flow or provide resistance. For short pipes with an appreciable incline as well as long pipelines with a slight inclination, this effect could lead to significant changes in the expected pressures. 
Moreover, gravity could encourage the development of larger velocity gradients and increase the impact of the inertial effect.
Thus, in order to judge the significance of gravity and inertial effects, one needs to undertake the challenging task of performing a numerical simulation  of  the flow in a network with realistic topology and geometry.
This task is the subject of the present article. 

In order to explain why this is a challenging task,  we observe that in the absence of gravity and inertial effects under steady-state conditions, the ordinary differential equation (ODE) governing the spatial distribution of pressure/density  in a pipe can be integrated analytically to yield a closed-form algebraic expression for the first integral involving the pressures at the ends of the pipe as well as the (constant) flow through it, and over a network, the problem then reduces to a system of algebraic equations \cite{Dvijotham2015Jun,Singh2020, nrsolver, ss-soln-existence, network-partitioning} involving the pressures at the network junctions and the flows along the edges of the network. 
This significant mathematical simplification has converted the problem into the search for the zero of a  finite-dimensional nonlinear algebraic system.
Additionally, we show (in Section~\ref{sec:closed-form}) that as long as the ideal gas equation of state holds, a closed-form first integral can be derived analytically irrespective of the presence or absence of the terms corresponding to inertial or gravitational effects.
However, if one considers a non-ideal equation of state, then the inclusion of either the inertial or gravitational effect or both renders this task hopeless, for  we now have a system of first-order  ODEs (one for each pipe)  which are coupled through algebraic equations in the form of boundary conditions on the pressure or a balance of mass flows at either end. 
One of our main contributions in this article is the demonstration of how it is yet possible to compute the solution of this nonlinear differential-algebraic system via the search for zeros of a finite-dimensional nonlinear system by use of the Newton-Raphson algorithm.

Our stratagem is to construct, in lieu of a putative first integral of each ODE, an expression involving the pressures at the ends of the pipe as well as the flow through it such that the expression behaves as a residual of the  solution. This expression is constructed using the abstract solution operator of the ODE. In the  use of the Newton-Raphson iterations, the Jacobian terms are then evaluated by solving the corresponding ODEs that govern forward sensitivity.
Moreover, at the outset, we propose a variable transformation that alleviates the poor scaling of the ODE, and we also introduce a two-point collocation scheme as a coarse approximation from which to find initial guesses for the Newton iterations. 
Subsequent simulation studies  were conducted for a single pipe, as well as a large-scale pipeline network with real data.

Through our work in this article, we would like to
(i) make a definitive evaluation of the importance of gravity and inertial effects in the flow of gas in pipeline networks, and delineate regimes of flow where these effects are negligible, and (ii) illustrate a computational framework for the simulation of such flows in large-scale networks when these effects are non-negligible and deserve consideration.

Although the presentation in this article is motivated by the steady-state flow of natural gas, the computational framework we develop is not specific to gas pipeline networks. The essential mathematical structure consists of balance laws  with nonlinear constitutive relations and gravitational effects formulated on the edges of a network which are coupled through algebraic constraints at junctions. Such formulations arise naturally in a variety of pipeline and thermal-fluid network models, including crude oil and refined product pipelines, hydrogen transport networks, CO\textsubscript{2} pipelines, and district heating systems, each of which exhibits its own equation of state, friction terms, and driving forces. In all such cases, the governing equations lead to analogous network-coupled boundary-value problems that resist reduction to purely algebraic models. The methodology proposed here provides a general and flexible approach for computing steady states of these systems without requiring algebraic reduction of governing equations, and is therefore applicable well beyond the specific setting  of natural gas that is considered in this work.

%% file: sections/pipe.tex
\section{Steady flow in a pipe} \label{sec:governing_equations}
Gases differ from liquids in that their densities can be changed appreciably by changing the pressure and/or temperature. The Equation of State (EoS) of a gas quantifies this relationship mathematically. Our interest is restricted to isothermal processes, i.e., one where temperature is constant, so the EoS relating pressure $p$ and density $\rho$ can be expressed in either of the forms
\begin{equation}
    p = p(\rho) \; \mathrm{or} \; \rho = \rho(p)
    \label{eqn:EOS}
\end{equation}
For an incompressible fluid like water, for instance, the EoS expresses the  fact that the fluid density $\rho$ is independent of pressure, i.e., a constant, while for gases, one commonly uses a non-ideal EoS (the CNGA EoS \cite{GyryaZlotnik2017})  of the form
    

\begin{flalign}
    \rho(p) = \frac{b_1 p + b_2 p^2}{R_gT}. \label{eq:cnga}
\end{flalign}
where, $b_1$ and $b_2$ are gas and temperature-dependent constants while $R_g$ is the applicable gas constant.
In the special case that $b_2 = 0, b_1 =1$, it is called the ideal gas EoS.
The values of $b_1$ and $b_2$ for the non-ideal (CNGA)  EoS  are given by the following expressions:
\begin{align}
  b_1 = 1 + \left( \dfrac{p_{atm}}{6894.75729}\right) 
  \left( \dfrac{a_1 10 ^ {a_2 G}}{(1.8T) ^ {a_3}} \right) ~~(\text{\si{unitless}}),\label{eq:b1}\\
    b_2 = \left(\dfrac{1}{6894.75729}\right)
    \left(\dfrac{a_1 10 ^ {a_2 G}}{(1.8T)^{a_3} } \right) ~~(\text{\si{\per\pascal}}). \label{eq:b2}
\end{align}
Here, $b_1$ and $b_2$ are calculated in terms of other non-dimensional constants $a_1 = 344400$, $a_2 = 1.785$, $a_3 = 3.825$, specific gravity of natural gas $G = 288.706$ and atmospheric pressure $p_{atm} = 101350\; \si{\pascal}$. 

The ``isothermal'' speed of sound, denoted by $c$, is related to the EoS through
\begin{equation}
    c(\rho) \triangleq \sqrt{p'(\rho)}
\end{equation}
For the ideal gas EoS, the speed of sound $c(\rho)$ is a constant.

Let the spatial domain coincide with the axis of a rectilinear pipe of length $L$ and constant cross-sectional area $A$, with origin chosen to be the left end, i.e., we assume 
the velocity depends on the axial coordinate alone.

Given our interest in the steady-state flow equations, the local acceleration (time derivative) terms are dropped and  
density $\rho(x)$, pressure $p(x)$, and velocity $v(x)$ are governed by the equations representing the  steady-state balance of mass and momentum for the flow of a gas in a pipe
which read
\begin{subequations}
\begin{gather}
  (\rho v)_x = 0 \quad \mathrm{on} \;  (0, L), \\ 
  \rho vv_x + p_x  = -\dfrac{\lambda}{2D}v|v| + \rho g e_{\parallel} \quad \mathrm{on} \; (0, L),
  \end{gather}
  \label{eqn:balance-original-NS}
\end{subequations}
where  $\lambda, D$ are positive constants representing the friction coefficient and diameter of the pipe, respectively, while $e_{\parallel} \in [-1, 1]$ is the component of the gravitational unit vector along the assumed flow direction of the pipe. Thus, for horizontal pipes, $e_{\parallel} =0$ while for vertical pipes $e_{\parallel} = \pm1$ depending on whether the assumed flow direction is upwards or downwards. If $\theta \in [-\pi/2, \pi/2]$ is the angle of inclination of the pipe with respect to the horizontal, then $e_{\parallel} = \sin \theta$ for an appropriate value of $\theta$. The system \eqref{eqn:balance-original-NS} can be rewritten  in terms of the constant mass flow $f = \rho v A$ as 
\begin{equation}
  \left(\dfrac{f^2}{\rho A^2}\right)_x + p_x  = -\dfrac{\lambda}{2DA^2}\dfrac{f |f|}{\rho} + \rho g \sin \theta\quad \mathrm{on} \; (0, L). 
  \label{eqn:momentum}
\end{equation}
Note that in the above equation,  the inertial term is the first one on the left.

\subsection{Nondimensional equations}
In order to render the governing equation \eqref{eqn:momentum} dimensionless we use nominal values (indicated by subscript 0) for the length, velocity, pressure and density to write
$$x = L_0 \bar{x}, \; v = v_0 \bar{v}, \; p = p_0 \bar{p}, \; \rho = \rho_0 \bar{\rho}.$$
We pick $L_0$ based on the length of the pipeline, $\rho_0$ can be set to its value at standard temperature and pressure, while $p_0, v_0$ are chosen based on the average flow conditions known to prevail in pipelines.
We set the nominal cross-sectional area $A_0=1$ so that the nominal mass flow $f_0 = \rho_0 v_0 A_0$, and we set $f = f_0 \bar{f}$.
We use $c_0 = c(\rho_0)$ as the nominal speed of sound, so that $c(\rho) = c_0\bar{c}(\bar{\rho})$. For the ideal gas EoS, $\bar{c} = 1$.

If we recognize the Mach number $\mathrm{M} = v_0/c_0$, the Euler number $\mathrm{Eu} = p_0/(\rho_0c_0^2)$ and the Froude number $\mathrm{Fr} = v_0/\sqrt{gL_0}$ and designate the nondimensional group $R_1 = \mathrm{M}^2/(\mathrm{Eu} \bar{A}^2)$, $R_2 = \mathrm{M}^2/(\mathrm{Eu}{\mathrm{Fr}^2})$, and $\beta=\lambda/(2\bar{D})$ then we obtain
\begin{equation}
  R_1\left(\dfrac{\bar{f}^2}{\bar{\rho}}\right)_{\bar{x}} + \bar{p}_{\bar{x}}  = -R_1\beta \dfrac{\bar{f} |\bar{f}|}{\bar{\rho}} + R_2\bar{\rho} \sin \theta\quad \mathrm{on} \; (0, L/L_0). 
  \label{eqn:nondim-momentum}
\end{equation}
The EoS, $\rho = \rho(p)$, is transformed to $\bar{\rho} = \bar{\rho}(\bar{p})$.
If we now drop the overbars in the notation for convenience, with the recognition that all quantities appearing henceforth are nondimensional,  we can use the EoS to express the spatial derivative of density and obtain 
\begin{equation}
    p_x = \rho\dfrac{\left( R_2\rho^2\sin \theta - R_1 \beta f |f| \right)}{ \left( \rho^2 -R_1f^2\rho'(p)\right)}.
    \label{eqn:p-ode-nondim}
\end{equation}
For subsequent use, we will denote 
\begin{equation}
    G(p, f) \triangleq \rho\dfrac{\left( R_2\rho^2\sin \theta - R_1 \beta f |f| \right)}{ \left( \rho^2 -R_1f^2\rho'(p)\right)}
    \label{eqn:G}
\end{equation}
\subsection{Closed-form first integrals for ideal gas EoS}
\label{sec:closed-form}
For a pipe with governing equation given by the ODE \eqref{eqn:p-ode-nondim}, two sets of boundary conditions are mathematically sensible.
Without loss of generality, the initial condition $p(0) = p_0$ could be specified at one end of the pipe while at the other end, either the injection or withdrawal is specified (i.e., $f$ given) or $p(L) = p_L$ is known.
In the former case, the ODE is now fully specified with $f$ given, and if a solution exists, the ODE can be solved to find $p(x)$ for $x \in (0, L]$. In the latter case, a procedure analogous to the  ``shooting method'' can be employed to find the value of $f$---if it exists---which, in conjunction with the initial condition $p(0) = p_0$, yields the given value at $p(L)$. Thereafter, $p(x)$ can be determined by solving the ODE with $f$ and $p_0$.

In both cases, however, it is useful to regard $p = p(x, p_0, f)$, recognizing the dependence of the solution on the initial condition as well as the parameter.
While  a closed-form first integral of \eqref{eqn:p-ode-nondim} cannot be derived  for a general EoS, such a derivation is possible for the case of the ideal gas EoS, wherein 
$$\rho(p) = \mathrm{Eu}\cdot p\; \textrm{while}\; \rho'(p) = \mathrm{Eu}.$$
In that case,  defining 
$$\hat R_1 = \frac{R_1}{\mathrm{Eu}},\quad \hat R_2 = R_2 \mathrm{Eu}, $$ 
we can enumerate four cases as follows, depending on whether the model includes the terms due to inertia and gravity: \\

\noindent \textbf{No gravity and no inertia} -- The ODE in this case is 
\begin{gather*}
p_x = -R_1 \beta \frac{f|f|}{\rho}
\end{gather*}
whose integration yields 
\begin{equation}
    p^2_0 - p^2_L  -2L\hat R_1 \beta f|f| = 0.
            \label{eqn:no-gravity-no-inertia}
\end{equation}

\noindent \textbf{With inertia only} -- The ODE is now $$p_x = \dfrac{\left( - \rho R_1 \beta f |f| \right)}{ \left( \rho^2 -R_1f^2\rho'(p)\right)}$$ which has the first integral
\begin{equation}
p^2_0 - p^2_L - \hat R_1 f^2 \ln{\frac{p^2_0}{p^2_L}}  -2L\hat R_1 \beta f|f| = 0.
\label{eqn:inertia-only}
\end{equation}  

\noindent \textbf{With gravity only} -- The ODE   $$p_x =  \dfrac{\left( R_2\rho^2\sin \theta - R_1 \beta f |f| \right)}{ \rho}$$ integrates to 
    \begin{equation}
        e^{\gamma}p^2_0 - p^2_L  -2L\hat R_1 \beta f|f| \left( \dfrac{e^{\gamma} - 1}{\gamma} \right) = 0, \; 
            \label{eqn:gravity-only}
    \end{equation}
    where, $\gamma \triangleq  2L\hat R_2 \sin \theta$. \\
    
\noindent \textbf{With gravity and inertia} -- The ODE is \eqref{eqn:p-ode-nondim}, and using $$\int\dfrac{(y^2 -a)}{y(y^2-b)}dy = \dfrac{(b-a)\log{(y^2-b)} + 2a\log{y}}{2b},$$
    the first integral has the form
    \begin{equation}
    \begin{gathered}
        \left(\hat R_1 f^2 - \delta \right)\log{ \dfrac{\left( p^2_0 - \delta \right)  }{\left( p^2_L - \delta \right)}} - 
         \hat R_1 f^2 \log{\frac{p^2_0}{p^2_L}}  -2L\hat R_1 \beta f|f| = 0 
    \end{gathered}
    \label{eqn:with-both}
    \end{equation}
    where, $\delta = \beta \hat R_1 f|f|/(\hat R_2 \sin \theta)$.

 An observation about these integrated forms is that in the limit of zero inertia for example, \eqref{eqn:with-both} yields \eqref{eqn:gravity-only} while \eqref{eqn:inertia-only} reduces to \eqref{eqn:no-gravity-no-inertia} and similarly, for zero gravity, \eqref{eqn:gravity-only} reduces to \eqref{eqn:no-gravity-no-inertia}.

%% file: sections/network.tex
\section{Steady flow in a network of pipelines}
We have seen that the model for flow in a pipe takes the form of a first-order ODE in general, unless its first integral can be derived in closed form. However, the subject of this article is flow in a network of pipelines, to which we now turn. In practice,
it is not \emph{pipes} alone that constitute a network, but \emph{compressors} are also necessary to boost the pressure and compensate for frictional losses and a decrease in pressure as the fluid flows long distances. For the purposes of simulation, it is sufficient to model both pipes and compressors as one-dimensional objects since we have already seen in Section~\ref{sec:governing_equations} that flow in a pipe is assumed
to depend on the axial coordinate alone. Thus, the network formed by the pipes and compressors can be thought of as a graph, with the edges of the graph being the pipes/compressors and the nodes being intersections of two or more edges.

Accordingly, let $G = (V, E)$ denote the graph of the network, where $V$ and $E$ are the corresponding sets of \emph{nodes} and \emph{edges} of the network. The symbols $|V|$ and $|E|$ will denote the cardinality of the sets $V$ and $E$, respectively.
The edge set is given by $E = P \bigcup C$ where $P\neq \emptyset$ and $C$ are the set of \emph{pipe} and \emph{compressor} elements, respectively.  
Thus,  nodes $i$ and $j$ are connected by an element $(i, j)$ that belongs to either $P$ or $C$. 

For each node $i \in V$, we let $p_i \in \mathbb{R}$ and $\rho_i \in \mathbb{R}$ denote the nodal value of the pressure and density, respectively while $q_i \in \mathbb{R}$ is the nodal injection or withdrawal depending on whether it is positive or negative.

For each edge $(i,j) \in P\subseteq E$, we associate positive constants representing its cross-sectional area $A_{ij}$, diameter $D_{ij}$, length $L_{ij}$, and friction factor $\lambda_{ij}$. In addition, we have the functions representing pressure and density $p^{ij}(x), \rho^{ij}(x)$  defined over the domain $[0, L_{ij}]$ as well as the constant mass-flow $f^{ij}$. 
The nodal quantities (pressure/density in this case) are coupled to the quantities defined along each edge via the mathematical assumption of continuity. While this is the most common coupling condition used, other coupling conditions are possible, as outlined in \cite{Herty-coupling}.

However, for an edge $(i,j) \in C \subset E$, we have an associated flow $f^{ij} \in \mathbb{R}$ and a compressor-ratio $\alpha_{ij} \geqslant 1$. In other words, a compressor is assumed to act by giving a multiplicative boost to the pressure through its specified compressor-ratio. 

Each vertex in $V$ is either a slack vertex (density/pressure specified) or a non-slack vertex (injection given). In other words $V = V_s \bigcup V_{ns}$.
At a vertex, the algebraic sum of the flow along every incoming and outgoing edge will equal the injection (known or unknown) at the given vertex. 

Therefore, the  system of differential-algebraic network flow equations in this case  consists of the $|V| + |E| + 2|P|$ governing equations and the same number of unknowns for flow along each edge (pipe/compressor) and coupling conditions for each junction  as follows ($*$ superscript indicates given data)
\begin{flalign}
~~\forall (i,j) \in C\subset E, \quad
\alpha^{*}_{ij}p_i - p_j = 0, \label{eqn:NF-C}
\end{flalign}
\begin{subnumcases}{\label{eqn:NF-P} \forall (i,j) \in P\subseteq E:} 
    p^{ij}_x = G(p^{ij}, f^{ij}) \; \mathrm{on}\; (0, L^*_{ij}], \label{eqn:p-ode}\\
    p^{ij}(0) = p_i, \label{eqn:p-ic}\\
    p^{ij}(L^*_{ij}) = p_j,
\end{subnumcases}
%
\begin{subequations}
\label{eqn:NF-bc}
    \begin{flalign}
       & ~~\forall j \in V_s, \quad p_j = p^*_j,  \\
       & ~~\forall j \in V_{ns}, \; \sum_{(i, j) \in E} f^{ij} - \sum_{(j, i) \in E} f^{ji}  = q^*_j.
    \end{flalign}
\end{subequations}
\begin{definition}{(Network system)}
    \label{def:network-syst}
The network system associated with graph $G(V,E)$ is the differential-algebraic system of equations  \eqref{eqn:NF-C}, \eqref{eqn:NF-P}, and \eqref{eqn:NF-bc}.
\end{definition}

If the  network system (Definition~\ref{def:network-syst})  is solved, i.e., $p^{ij}(x)$ and $f^{ij}$ are known for every $(i,j) \in P \subset E$ and $(i,j) \in E$ respectively, while $p_i$ is known for every $i \in V_{ns}$, then it is possible to evaluate $q_j$ for $j\in V_s$---the counterpart of the known quantity $p^*_j$ in \eqref{eqn:NF-bc}. 
\begin{remark}
\label{rem:network-syst}
The preceding statement yields the insight that the crux of the problem in the solution of the network system in Definition~\ref{def:network-syst} is the determination of the junction pressures and the edge flows that ensure consistency of the ODEs \eqref{eqn:NF-P} defined on the edges as well as the other equations \eqref{eqn:NF-C} and \eqref{eqn:NF-bc}.
\end{remark}


In the preceding section, we already saw in equations \eqref{eqn:no-gravity-no-inertia}, \eqref{eqn:inertia-only}, \eqref{eqn:gravity-only}, \eqref{eqn:with-both} that for the case of the ideal gas EoS, \eqref{eqn:NF-P} can be integrated analytically to yield a closed-form first integral in the form $F(p_i, p_j, f^{ij}) = 0$. In such a case, the complete network system (Definition~\ref{def:network-syst}) can be represented equivalently as a system of nonlinear equations as follows:
\begin{subequations}
    \begin{align}
        &  \alpha^{*}_{ij}p_i - p_j = 0 \quad \forall (i,j) \in C, \\
        & F(p_i, p_j, f^{ij}) = 0 \quad \forall (i,j) \in P, \\
    & \sum_{(i, j) \in E} f^{ij} - \sum_{(j, i) \in E} f^{ji}  = q^*_j \quad \forall j \in V_{ns},\\
    & p_j = p^*_j \quad \forall j \in V_{s}.
    \end{align}
    \label{eqn:integrated-network-syst}
\end{subequations}
Then, the problem is to determine the $|V| + |E|$ number of unknowns associated with the junctions and the edges as observed in Remark~\ref{rem:network-syst}.
The form of the problem most often considered is when $F(p_i, p_j, f^{ij})$ has the form in equation~\eqref{eqn:no-gravity-no-inertia}, and many of the cited works, such as \cite{Dvijotham2015Jun,ojha2017solving,Singh2019,Singh2020,nrsolver,ss-soln-existence,network-partitioning} are in that context.

However, in the case of a non-ideal EoS, these developments are to no avail; one cannot derive the first integral of the ODE  to convert the network system (Definition~\ref{def:network-syst}) into  \eqref{eqn:integrated-network-syst}. Instead, the system has to be recognized as a non-linear differential-algebraic system consisting of first-order ODEs defined on edges of the network coupled with the nodal equations in \eqref{eqn:NF-bc}, the values of the pressure at the node, and the mass flow on the edges.

\subsection{Solution approach}
\label{sec:problem-class}
The network system in Definition~\ref{def:network-syst} bears a resemblance to the canonical form of a semi-explicit order one differential-algebraic equation (DAE). However, on closer examination, it will be seen that the Jacobian of the algebraic part of the DAE is not invertible, thus violating the hypothesis required to bring to bear the established techniques for numerical solution of DAEs \cite{Gear1971Jan,Ascher1989Jun,Gear2006Jul}. Another point of difference is that even the initial conditions of the system in Definition~\ref{def:network-syst} are unknown and need to be determined as part of the solution, so that the problem defies classification into a known class within DAEs.

Our proposed solution approach instead draws inspiration from Remark~\ref{rem:network-syst} and \eqref{eqn:integrated-network-syst}, the case when the ODE has a closed-form first integral of the form $F(p_i, p_j, f^{ij})$ which acts as a residual of the solution. In that case, one uses the Newton-Raphson algorithm  to 
determine the $|V| + |E|$ unknowns through the Jacobian whose elements include the partial derivatives $$\doe{F}{p_i}, \doe{F}{p_j}, \doe{F}{f^{ij}}.$$

We posit that even though a closed-form first integral cannot be derived, the problem at hand may still be cast as the determination of the $|V| + |E|$ unknowns as in \eqref{eqn:integrated-network-syst} associated with the junctions and the edges if a suitable expression for a residual can be constructed, i.e., in place of the first integral we require a different expression $F(p_i, p_j, f^{ij})$ in \eqref{eqn:integrated-network-syst}.
We propose an expression that measures the mismatch at one end between the specified pressure and predicted pressure from the ODE solution as a function of the flow and the specified pressure (initial condition) at the other end.
Specifically, denoting the solution of the ODE \eqref{eqn:p-ode} with the initial condition \eqref{eqn:p-ic} as $p^{ij}(x, p_i, f^{ij})$, we recognize  that $p^{ij}(L^*_{ij}, p_i, f^{ij}) = p_j$ at a solution.

In order to  construct a residual, we choose a smooth algebraic function  $R(x_1, x_2)$ with the property that $$R(x_1, x_2) = 0 \implies x_1 = x_2.$$ 
Some obvious choices for $R(x_1, x_2)$ are $x_1 - x_2, \; (x_1 - x_2)^2, \; x_1^3 - x_2^3$ etc.
Thus, the equation $$R(p^{ij}(L^*_{ij}, p_i, f^{ij}), p_j) = 0  \quad \equiv \quad p^{ij}(L^*_{ij}, p_i, f^{ij}) = p_j.$$
Hence, continuing to use the notation established for the closed-form first integral,  we set
\begin{equation}
    F(p_i, p_j, f^{ij}) = R(p^{ij}(L^*_{ij}, p_i, f^{ij}), p_j)
    \label{eqn:ode-residual}
\end{equation}
 Our strategy is to  solve  \eqref{eqn:integrated-network-syst} iteratively with the Newton-Raphson
algorithm as before using \eqref{eqn:ode-residual}.

In order to do so, the Jacobian terms corresponding to $F(p_i, p_j, f^{ij})$ in  \eqref{eqn:ode-residual} can be computed as 
\begin{subequations}
\begin{gather}
    \doe{F}{p_i} = s^{ij}_p(L^*_{ij})\doe{R}{x_1}\left(x_1 = p^{ij}(L^*_{ij}, p_i, f^{ij}),p_j\right), \\
    \doe{F}{f^{ij}} = s^{ij}_f(L^*_{ij})\doe{R}{x_1}\left(x_1 = p^{ij}(L^*_{ij}, p_i, f^{ij}),p_j\right),\\
    \doe{F}{p_j} =  \doe{R}{x_2}\left(p^{ij}(L^*_{ij}, p_i, f^{ij}), x_2 =p_j\right),
\end{gather}
\end{subequations}
where
    \begin{align}
        & s^{ij}_p(x) = \doe{p^{ij}(x, p_i, f^{ij})}{p_i},\;
        s^{ij}_f(x) = \doe{p^{ij}(x, p_i, f^{ij})}{f^{ij}}
    \end{align}
are the sensitivities of $p^{ij}$ with respect to the initial condition and the flow.

Denoting  the partial derivatives of $G$ with respect to its first and second arguments as $G_p, G_f$ respectively, the sensitivities are themselves governed by the forward-sensitivity system \cite{Dickinson1976Jun}:
\begin{subequations}
\begin{flalign}
    \frac{dp^{ij}}{dx}   & = G(p^{ij}, f^{ij}), \quad  p^{ij}(0) = p_i,\\
    \frac{ds^{ij}_p}{dx} & = G_p(p^{ij}, f^{ij}) s^{ij}_p, \quad  s^{ij}_p(0) = 1,\\
    \frac{ds^{ij}_f}{dx} & = G_p(p^{ij}, f^{ij}) s^{ij}_f + G_f(p^{ij}, f^{ij}), \quad s^{ij}_f(0) = 0.
\end{flalign}
\label{eqn:fwd-sensitivity}
\end{subequations}
When the system \eqref{eqn:integrated-network-syst} is now solved simultaneously using an iterative algorithm such as Newton-Raphson, derivative information in the form of the Jacobian comes with the added complexity that evaluation of the residual as well as the Jacobian in  \eqref{eqn:ode-residual} involves solution of an ODE system; \eqref{eqn:p-ode}, \eqref{eqn:p-ic} in the former and \eqref{eqn:fwd-sensitivity} in the latter.

\subsection{Reformulation for networks}
 In the context of steady-state solution for a network, the pressures and the flow in the course of iterations often traverse regions that render the ODE \eqref{eqn:p-ode}  stiff. Hence, it is advantageous to reformulate the ODE in a form conducive to obtaining steady-state pressures and flows for a  network of pipes and compressors.  
 To that end, consider an invertible and differentiable transformation 
 \begin{equation}
     p = \psi(\pi), \;\; \pi = \psi^{-1}(p)
     \label{eqn:diffeonorphism}
 \end{equation}
 that allows us to transform the ODE \eqref{eqn:p-ode} into
\begin{subequations}
    \begin{align}
     & \pi'(x) = H(\pi, f), \\
     & H(\pi, f) \triangleq \dfrac{G(\psi(\pi), f)}{\psi'(\pi)},\\
     & \pi(0) = \pi_0 = \psi^{-1}(p(0)).
 \end{align}
 \label{eqn:transformed-ODE}
\end{subequations}
 Restricting attention to monomials,  we found the transformation $\psi^{-1} : x\mapsto x^{3}$ that yields $\pi = p^3, p = \pi^{\frac{1}{3}}$ convenient  to solve the problem since it also allows the equation~\eqref{eqn:NF-C} to be recast in the form 
 $$\forall (i,j) \in C\subset E(G), \quad (\alpha^*_{ij})^3 \pi_i - \pi_j = 0.$$
 In this context, we note that even powers were a problem since the values of $\pi$ could become negative during the iterative process. 
 Thus, the system \eqref{eqn:integrated-network-syst} is now transformed to one where we now solve for the transformed variable $\pi_i$ at every node instead of $p_i$.

\subsection{Two-point collocation}
We recall that success with a Newton-Raphson iterative algorithm often depends critically on a judicious initial guess, and this dependence becomes more acute in our problem since it involves ODEs in the computation of both the residual and the Jacobian.
If we can replace the ODE by an approximate discrete formulation, it will serve as a coarser version of the problem as well as allow us to solve it to obtain a good initial guess for the transformed version of the system~\eqref{eqn:integrated-network-syst}.

We replace the transformed version of \eqref{eqn:NF-P} with the following approximation :
\begin{equation}
    \forall (i,j) \in P, \quad \dfrac{\pi_i - \pi_j}{L^*_{ij}} + \dfrac{ H(\pi_i, f^{ij}) + H(\pi_j, f^{ij}) }{2} = 0.
    \label{eqn:2-pt}
\end{equation}

There are two ways to interpret the discretized equation~\eqref{eqn:2-pt}, but we are motivated by the second explanation that follows.
The first explanation is that equation~\eqref{eqn:2-pt} is a finite-difference discretization that uses a trapezoidal rule to integrate the ODE. 
Viewed from this perspective, the justification is not compelling.
The second way to interpret equation~\eqref{eqn:2-pt} is as a two-point collocation \cite{Iserles}. 
By a solution $\pi(x)$ to the ODE, we usually mean that $\pi(x)$ should satisfy the ODE for all $x$ in the edge/pipe. However, collocation schemes relax this requirement and demand that the solution satisfy the ODE at a finite set of specified points only. Since $x = 0$ and $x= L^*_{ij}$ are the nodes of the network, we demand that the ODE be satisfied at those two points, resulting in the equation~\eqref{eqn:2-pt}. Thus, this is a natural and well-motivated replacement of the original ODE and can serve as a coarse problem to generate an initial guess.
Note that the contribution to the Jacobian in a Newton-Raphson iteration from equation~\eqref{eqn:2-pt} is easily calculated. 

%% file: sections/results.tex
\section{Results}
At the end of Section~\ref{sec:introduction}, we had stated that our aim was two-fold: to evaluate the importance of the inertial and gravitational effects, and also demonstrate the use of the proposed computational framework in the case of large-scale networks.  We accomplish these goals through three case studies.
In the first case study, the validity of the formulation is checked  on some of the standard GasLib instances \cite{Schmidt2017}. 
The second  study  allows examination of the effect of gravity and inertia in isolation by studying flow in a single pipeline while in the third case, we consider a pipeline network that spans multiple states of the continental US passing through both plains and mountains.

\subsection{Case study 1}
As a first step to check the validity of the formulation and solution approach proposed in Section~\ref{sec:problem-class} for  \eqref{eqn:integrated-network-syst}, we consider the case of networks with horizontal pipes   carrying an ideal gas. 
The assumption of  horizontal pipes is forced upon us because information  on pipe inclination or node elevation is not available in the test cases considered.
The solutions obtained from the formulation and solution approach stated in Section~\ref{sec:problem-class} should lead to a small residual \emph{when substituted in  the closed-form first integrals \eqref{eqn:no-gravity-no-inertia} and \eqref{eqn:inertia-only}}.
Moreover, this problem will also act as a test for the two-point collocation approximation.

Accordingly, we consider three cases from the GasLib instances \cite{Schmidt2017} and the residuals are tabulated in Table~\ref{tab:first-integral-residual}.
The values in the table offer a convincing argument  that the collocation scheme supplies a good initial guess for the ODE formulation. The values also indicate that inertial effects are negligible in these cases. 
This is not a surprising result in itself since gravity was absent and the same conclusion was reached in \cite{Osiadacz1984Jan}. Moreover, if we examine  the terms in \eqref{eqn:inertia-only}, we see that the term $\hat R_1 f^2 \ln({p^2_0}/{p^2_L})$ cannot be large for $\hat R_1 \sim O(M^2)$ with $M\ll 1$ for typical pipeline operation.
However, in order to investigate if the inertial effects may be neglected altogether in steady-state computations, we will focus on flow in a single pipe while considering gravity.


\begin{table}[htb]
\small
\caption{The table shows the maximum error/residual (among all pipes)  of the solution to the collocation and ODE formulations when substituted into the first-integrals. From the values, we observe that the collocation scheme supplies a good initial guess for the ODE formulation. The values also indicate that inertial effects are negligible in this case.}
\label{tab:first-integral-residual}
\centering
\begin{filecontents*}{data/gaslib-errors.csv}
networks,no-inertia-(colloc.),no-inertia-(ode),with-inertia-(colloc.),with-inertia-(ode)
GasLib-11,8.2E-05,8.2E-06,8.2E-05,8.2E-06
GasLib-40,9.2E-03,5.1E-04,9.2E-03,5.1E-04
GasLib-134,6.9E-02,3.7E-03,6.9E-02,3.7E-03
\end{filecontents*}
\csvreader[
  tabular={@{}lcccc@{}},
  table head=\toprule
    \multirow{3}{*}{Networks} & \multicolumn{2}{c}{No inertia, Eqn \eqref{eqn:no-gravity-no-inertia}} & \multicolumn{2}{c}{With inertia, Eqn. \eqref{eqn:inertia-only}} \\ \cmidrule(lr){2-5} 
     & \multicolumn{2}{c}{Max. error} & \multicolumn{2}{c}{Max. error} \\ \cmidrule(lr){2-5}
    & \multicolumn{1}{l}{Colloc.} & ODE & \multicolumn{1}{l}{Colloc.} & ODE \\ \midrule,
  late after line=\\,
  table foot=\bottomrule
]{data/gaslib-errors.csv}{}%
{\csvcoli & \num{\csvcolii} & \num{\csvcoliii} & \num{\csvcoliv} & \num{\csvcolv}}
\end{table}

 \subsection{Case study 2}
\begin{figure}[htb]
    \centering
    \includegraphics[scale=0.7]{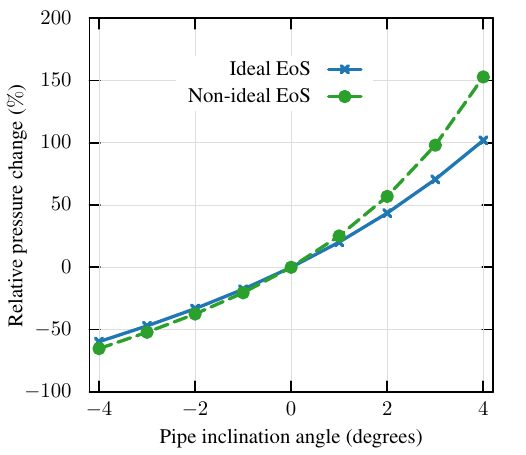}
    \caption{The figure shows the relative pressure change at the withdrawal end with respect to the pressure computed for a horizontal pipe. For both the ideal and non-ideal EoS, the gravitational effects are  significant and the relative change in pressure is a nonlinear function of the pipe inclination. We can also confirm that the choice of the EoS has  increasing  significance on the computational results as the (positive) angle of inclination increases}
    \label{fig:yamal}
\end{figure}

We now consider the flow in a single pipeline where pressure is specified at one end, and gas is withdrawn at a known rate at the other end in order to examine the effect of gravity and inertia thoroughly. 
Note that the pressure at the withdrawal end for a given value of withdrawal is dependent on the  length of the pipe, the inclination, as well as the pressure at the upstream end. However, by considering a single pipe, we are able to examine the variation of the pressure at the withdrawal end with respect to the pipe inclination for the same upstream pressure and pipe length. 
The problem parameters are chosen to correspond to those listed in  \cite{Brouwer2011Jun,osiadacz01} for the Yamal-Europe pipeline.
The pipeline is 122 kilometers long with a diameter of 1.422 meters, and a friction coefficient of 0.03. One end is maintained at a pressure of 8.8 MPa while gas is withdrawn at the other end at the rate of 400 kg/s.
We compute the pressure at the withdrawal end  for a variety of inclination angles, with and without the inertial effects, and for the ideal as well as a non-ideal EoS (see Fig.~\ref{fig:yamal}).
The inclination angles are varied around zero to ensure flow with and against gravity is considered. We choose the maximum inclination angle of 4 degrees from recognition of the fact that the highest elevation possible on the earth's surface is less than 9 kilometers.

From Table~\ref{tab:intertia-effects}, it is clear that the effects of inertia are negligible for each angle of inclination  for both the ideal and non-ideal EoS and due to this conclusion, subsequent results make no mention of whether inertial effects are incorporated or neglected. 

Figure~\ref{fig:yamal} examines the relative pressure change at the withdrawal end with respect to the pressure computed for a horizontal pipe for both the ideal and non-ideal EoS.
The gravitational effects are  significant and the relative change in pressure is a nonlinear function of the pipe inclination. 
Moreover, The figure also confirms that the choice of the EoS has  increasing  significance on the computational results  for larger  (positive) angles of inclination.

\begin{table}[htb]
\caption{Relative percentage change in pressure at the withdrawal end corresponding to the inclusion/exclusion of inertial effects for different angles of inclination.}
\label{tab:intertia-effects}
\centering
\begin{filecontents*}{data/inertia-errors.csv}
inclination,pressure-error-ideal,linepack-error-ideal,pressure-error-cnga,linepack-error-cnga
-4,7.67E-02,1.92E-02,9.79E-02,2.15E-02
-3,3.65E-02,1.11E-02,4.05E-02,1.18E-02
-2,1.83E-02,6.52E-03,1.86E-02,6.82E-03
-1,8.85E-03,3.58E-03,8.05E-03,3.61E-03
0,3.40E-03,1.54E-03,2.37E-03,1.27E-03
1,8.00E-05,3.99E-05,9.69E-04,6.14E-04
2,2.03E-03,1.10E-03,3.07E-03,2.28E-03
3,3.41E-03,2.00E-03,4.54E-03,3.90E-03
4,4.33E-03,2.71E-03,5.61E-03,5.53E-03
\end{filecontents*}
\csvreader[
  tabular={@{}ccc@{}},
  table head=\toprule
   Angle (deg.) & \multicolumn{1}{c}{Ideal EoS (\%)} & \multicolumn{1}{c}{Non-ideal EoS (\%)} \\ \midrule,
  late after line=\\,
  table foot=\bottomrule
]{data/inertia-errors.csv}{}%
{\csvcoli & \num{\csvcolii} & 
    \num{\csvcoliv} 
    }
\end{table}

\subsection{Case study 3}
\begin{figure}[htb]
    \centering
    \includegraphics[frame,width=0.4\linewidth]{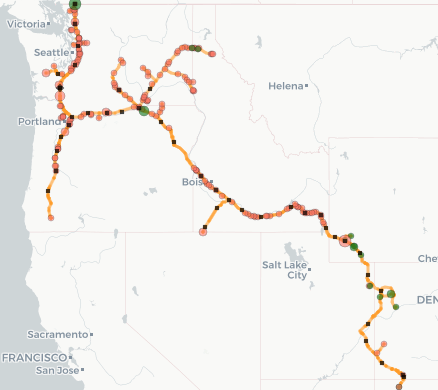}
    \caption{Northwest pipeline network. The pipeline is shown in orange, the black squares are the compressors, the red/green circles are the withdrawals/injections. The radius of the circle corresponds to the magnitude of injection/withdrawal. The network has a total length of 9228 km, comprising 977 nodes (of which 183  are injection/withdrawal nodes), 1033 edges, and 41 compressors.}
    \label{fig:NW}
\end{figure}

In this study, we shall  consider a large-scale network in the United States, the so-called Northwest pipeline network \cite{NW} that spans the states of Colorado, Utah, Idaho, Oregon, New Mexico, Wyoming, Montana, and Washington  as shown in Figure~\ref{fig:NW}. 
Unlike the first case study,  information on the elevation at given node locations is available, which  allows us to consider the effect of gravity. 
In addition,  it serves as a challenging and suitable test for our computational framework (see network metrics in Figure~\ref{fig:NW}).

In a network with given compressor ratios, the pressures at non-slack nodes and the flows in edges are determined by the slack pressure(s) as well as the non-slack injections/withdrawals. 
Thus one can interpret the effect of gravity only in a broad sense, unlike our careful comparison for a single pipe earlier. 
 
\begin{figure}[htb]
    \centering
    \includegraphics[scale=0.8]{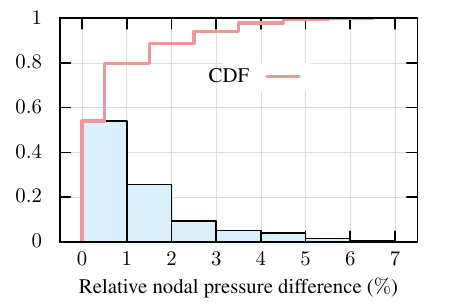}
    \caption{The abscissa shows the relative difference in predicted nodal pressures  due to gravity in the form of a (normalized) histogram as well as the empirical cumulative density function (CDF). The ideal gas EoS is assumed.}
    \label{fig:nw-ideal}
\end{figure}

\begin{figure}[htb]
    \centering
    \includegraphics[scale=0.8]{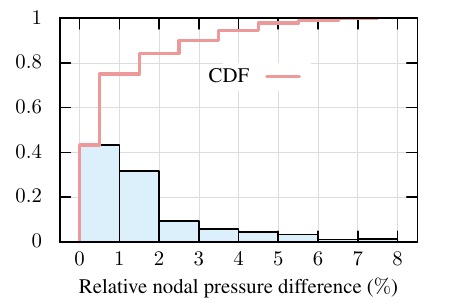}
    \caption{The abscissa shows the relative difference in predicted nodal pressures  due to gravity in the form of a (normalized) histogram as well as the empirical cumulative density function (CDF). The non-ideal  EoS is assumed.}
    \label{fig:nw-cnga}
\end{figure}

If we regard this problem from the viewpoint of a network operator, pipelines are usually operated close to their capacity limits for economic reasons, so a simulation is often used to  test if a given set of inputs could lead to pressures that violate safe operating regimes. Since the elevation information for nodes (or equivalently, pipeline inclination) is often missing, the  relative difference in predicted nodal pressures  that results from the inclusion of gravitational effects is a useful metric to examine.
These differences are summarized in the (normalized) histograms and cumulative density function (CDF) depicted in Figure~\ref{fig:nw-ideal} and Figure~\ref{fig:nw-cnga} for the ideal and non-ideal EoS respectively.

For the given network, it shows that at roughly  80\% of the nodes, the inclusion of gravity leads to a  relative pressure difference of less  than 2\%. 
Even though relative pressure differences of  7-8\% occur only at a handful of nodes, the operation of pipeline networks close to their capacity limits validates the need for higher-fidelity simulations that incorporate the gravitational effects.

%% file: sections/conclusion.tex
\section{Conclusion}
\label{sec:conclusion}
In this article, we have made a definitive evaluation of the importance of gravity and inertial effects in the steady flow of gas in pipeline networks, and delineated regimes of flow where these effects are negligible, and also illustrated a computational framework for the simulation of such flows in large-scale networks when these effects are non-negligible and deserve consideration.

We find that for the typical steady flow regimes in pipeline networks, irrespective of the presence or absence of gravitational effects, the inertial effects due to convective acceleration are insignificant and hence the corresponding term may be neglected in the flow equations.
On the other hand, the effects of gravity are significant in long pipes even with a few degrees of inclination.
We demonstrated a method to solve the nonlinear differential-algebraic system that arose, consisting of first-order ODEs (one for each pipe), which were coupled through algebraic equations in the form of boundary conditions on the pressure or a balance of mass flows at either end. 
We also introduced a two-point collocation scheme as a coarse approximation from which to find initial guesses for initiating the Newton-Raphson iterations.
The proposed computational framework was  successfully utilized in the simulation of steady flow for a large-scale pipeline network that spans multiple states in the continental US.

From a computational standpoint, a central outcome of this work is the demonstration that for flow models defined over networks, steady-state solutions  can be computed reliably even when the differential equation governing flow on the network edge cannot be simplified analytically to an algebraic relation. By defining residuals via solution operators and evaluating Jacobian information via forward sensitivity equations, Newton's method can be applied in a way that preserves the continuous structure of the governing equations while remaining scalable to large networks. The variable transformation and two-point collocation strategy further highlight how careful numerical design can mitigate stiffness and sensitivity to initial guesses in such problems. Taken together, these elements provide a practical template for solving steady-state boundary-value problems formulated on networks, where nonlinear physical laws  and geometric effects preclude closed-form equations on each edge.